\documentclass[a4paper,11pt]{article}
\usepackage{amssymb}
\usepackage{amsmath}
\usepackage[thinlines]{easybmat}

\renewcommand{\d}[1]{\;\mathrm{d}#1}

\newtheorem{lemma}{Lemma}

\renewcommand{\d}[1]{\;\mathrm{d}#1}

\title{CPHD filter derivation for extended targets}
\author{Umut Orguner}
\date{}
\begin{document}
\setlength{\interdisplaylinepenalty}{2500} \maketitle
\tableofcontents
\newpage
In this
document, an attempt to derive a cardinalized probability hypothesis
density filter (CPHD) (corrector formula only) for extended targets is given.
\section{Derivation}
We first
give some quantities related to the derivation.
\begin{itemize}
\item The extended target measurements are distributed according to
an iid cluster process. The corresponding likelihood is given as
\begin{align}
f(Z|x)=n!P_z(n|x)\prod_{z\in Z}p_z(z|x)
\end{align}
\item The false alarms are distributed according to an iid cluster
process also.
\begin{align}
f(Z_{FA})=n!P_{FA}(n)\prod_{z\in Z_{FA}} p_{FA}(z)
\end{align}
\item The multitarget prior $f(X_k|Z_{0:k-1})$ is assumed to be an iid cluster
process.
\begin{align}
f(X_k|Z_{0:k-1})=n!P_{k+1|k}(n)\prod_{x_k\in X_k} p_{k+1|k}(x_k)
\end{align}
where
\begin{align}
p_{k+1|k}(x_k)\triangleq N_{k|k-1}^{-1}D_{k|k-1}(x_k)
\end{align}
with $N_{k|k-1}\triangleq \int D_{k|k-1}(x_k)\d x_k$.
\end{itemize}
The p.g.fl corresponding to the updated multitarget density
$f(X_k|Z_k)$ is then given as
\begin{align}
G_{k|k}[h]=\frac{\frac{\delta}{\delta Z_k}F[0,h]}{\frac{\delta}{\delta Z_k}F[0,1]}
\end{align}
where
\begin{align}
F[g,h]\triangleq \int& h^{X}G[g|X]f(X|Z_{k-1})\delta X\\
G[g|X]\triangleq \int& g^Zf(Z|X)\delta Z
\end{align}
\begin{itemize}
\item Calculation of $G[g|X]$: Suppose given the target states $X$
the measurement sets corresponding to different targets are
independent. Assume that targets are detected with probabilities
$p_D(x)$. Then the p.g.fl for a single targets measurements
become
\begin{align}
G[g|x]=1-p_D(x)+p_D(x)G_Z[g|x]
\end{align}
Then p.g.fl for the measurements belonging to all targets become
\begin{align}
G[g|X]=(1-p_D(x)+p_D(x)G_Z[g|x])^X
\end{align}
With the addition of false alarms we have
\begin{align}
G[g|X]=G_{FA}[g](1-p_D(x)+p_D(x)G_Z[g|x])^X
\end{align}
\item Calculation of $F[g,h]$: substituting $G[g|X]$ above into the
definition of $F[g,h]$ we get
\begin{align}
F[g,h]=&\int h^XG_{FA}[g] (1-p_D(x)+p_D(x)G_Z[g|x])^Xf(X|Z_{k-1})\delta X\\
=&G_{FA}[g]\int  \left(h(1-p_D(x)+p_D(x)G_Z[g|x])\right)^Xf(X|Z_{k-1})\delta X\\
=&G_{FA}[g]G_{k|k-1}\big[h(1-p_D(x)+p_D(x)G_Z[g|x]\big]\\
=&G_{FA}[g]G_{k|k-1}\big[h(1-p_D+p_DG_Z[g]\big]
\end{align}
where we omitted the arguments $x$ of the functions.
\end{itemize}

For the iid cluster processes, we know the following identities
\begin{align}
G_{FA}[g]=&G_{FA}(p_{FA}[g])\\
G_{Z}[g|x]=&G_{Z}(p_{z}[g]|x)
\end{align}
where the functions $G_{FA}(\cdot)$ and $G_Z(\cdot|x)$ are the
probability generating functions for the cardinality distributions
$P_{FA}(\cdot)$ and $P_z(\cdot|x)$ respectively. We also know the
following derivative expressions.
\begin{align}
\frac{\delta}{\delta Z}G[g]=G^{(n)}(p[g])\prod_{z\in Z}p(z)
\end{align}
Now taking the derivatives of $F[g,h]$ with respect to $Z$.
\begin{align}
\frac{\delta}{\delta z_1}F[g,h]=&G_{FA}^{(1)}(p_{FA}[g])p_{FA}(z_1)G_{k|k-1}\Big(p_{k|k-1}\big[h(1-p_D+p_DG_Z(p_z[g]))\big]\Big)\nonumber\\
&+G_{FA}(p_{FA}[g]) G_{k|k-1}^{(1)}\Big(p_{k|k-1}\big[h(1-p_D+p_DG_Z(p_z[g]))\big]\Big)p_{k|k-1}\big[hp_DG_Z^{(1)}(p_z[g])p_z(z_1)\big]\nonumber\\
=&F[g,h]p_{FA}(z_1)\bigg( \frac{G_{FA}^{(1)}(p_{FA}[g])}{G_{FA}(p_{FA}[g])}\nonumber\\
&+\frac{G_{k|k-1}^{(1)}\Big(p_{k|k-1}\big[h(1-p_D+p_DG_Z(p_z[g]))\big]\Big)}{G_{k|k-1}\Big(p_{k|k-1}\big[h(1-p_D+p_DG_Z(p_z[g]))\big]\Big)}p_{k|k-1}\Big[hp_DG_Z^{(1)}(p_z[g])\frac{p_z(z_1)}{p_{FA}(z_1)}\Big]\bigg)\nonumber\\
=&F[g,h]p_{FA}(z_1)\bigg(\zeta_{FA}^{(1)}(p_{FA}[g])\nonumber\\
&+\zeta_{k|k-1}^{(1)}\Big(p_{k|k-1}\big[h(1-p_D+p_DG_Z(p_z[g]))\big]\Big)p_{k|k-1}\Big[hp_DG_Z^{(1)}(p_z[g])\frac{p_z(z_1)}{p_{FA}(z_1)}\Big]\bigg)\label{Eqn:FirstDerivativeF}
\end{align}
where we use the functions $\zeta^{(i)}_{FA}(\cdot)$ and
$\zeta_{k|k-1}^{(i)}(\cdot)$ defined as
\begin{align}
\zeta^{(i)}_{FA}(x)\triangleq&\left(\frac{G_{FA}^{(1)}(x)}{G_{FA}(x)}\right)^{(i-1)}=\log^{(i)} G_{FA}(x)\\
\quad\zeta_{k|k-1}^{(i)}(x)\triangleq&\left(\frac{G_{k|k-1}^{(1)}(x)}{G_{k|k-1}(x))}\right)^{(i-1)}=\log^{(i)} G_{k|k-1}(x)
\end{align}
The second order derivative is given as
\begin{align}
\frac{\delta}{\delta z_1\delta z_2}F[g,h]=&F[g,h]\bigg(\prod_{z'=z_1,z_2}p_{FA}(z')\bigg)\Bigg(\prod_{z'=z_1,z_2}\bigg(\zeta_{FA}^{(1)}(p_{FA}[g])\nonumber\\
&+\zeta_{k|k-1}^{(1)}\Big(p_{k|k-1}\big[h(1-p_D+p_DG_Z(p_z[g]))\big]\Big)p_{k|k-1}\Big[hp_DG_Z^{(1)}(p_z[g])\frac{p_z(z')}{p_{FA}(z')}\Big]\bigg)\nonumber\\
&+\zeta_{FA}^{(2)}(p_{FA}[g])+\zeta_{k|k-1}^{(2)}\Big(p_{k|k-1}\big[h(1-p_D+p_DG_Z(p_z[g]))\big]\Big)\nonumber\\
&\times\prod_{z'=z_1,z_2} p_{k|k-1}\Big[hp_DG_Z^{(1)}(p_z[g])\frac{p_z(z')}{p_{FA}(z')}\Big]\nonumber\\
&+\zeta_{k|k-1}^{(1)}\Big(p_{k|k-1}\big[h(1-p_D+p_DG_Z(p_z[g]))\big]\Big)p_{k|k-1}\Big[hp_DG_Z^{(2)}(p_z[g])\prod_{z'=z_1,z_2}\frac{p_z(z')}{p_{FA}(z')}\Big]\Bigg)\label{Eqn:SecondDerivativeF}
\end{align}
Then we can write the general formula as
\begin{align}
\frac{\delta}{\delta Z}F[g,h]=&F[g,h]\bigg(\prod_{z'\in Z}p_{FA}(z')\bigg)\sum_{\mathcal{P}\angle Z}\prod_{W\in\mathcal{P}}\bigg(\zeta_{FA}^{(|W|)}(p_{FA}[g])\nonumber\\
&+\sum_{\mathcal{Q}\angle W}\zeta_{k|k-1}^{(|\mathcal{Q}|)}\Big(p_{k|k-1}\big[h(1-p_D+p_DG_Z(p_z[g]))\big]\Big)\prod_{V\in \mathcal{Q}}\eta_V[g,h]  \bigg)\label{Eqn:FormulaForF}
\end{align}
where
\begin{align}
\eta_V[g,h]\triangleq p_{k|k-1}\Big[hp_DG_Z^{(|V|)}(p_z[g])\prod_{z'\in V}\frac{p_z(z')}{p_{FA}(z')}\Big]
\end{align}
A proof of \eqref{Eqn:FormulaForF} can be found in Section~\ref{Sec:ProofOfDerivative}.
Setting $g=0$ and taking derivative with respect to $x$.
\begin{align}
\frac{\delta}{\delta x}\frac{\delta}{\delta Z}F[0,h]=&\frac{\delta}{\delta Z}F[0,h]\zeta_{k|k-1}^{(1)}\Big(p_{k|k-1}\big[h(1-p_D+p_DG_Z(0))\big]\Big)\nonumber\\
&\times\big(1-p_D(x)+p_D(x)G_Z(0)\big)p_{k|k-1}(x)\nonumber\\
&+F[0,h]\bigg(\prod_{z'\in Z}p_{FA}(z')\bigg)\sum_{\mathcal{P}\angle Z}\Bigg(\prod_{W\in\mathcal{P}}\bigg(\zeta_{FA}^{(|W|)}(0)\nonumber\\
&+\sum_{\mathcal{Q}\angle W}\zeta_{k|k-1}^{(|\mathcal{Q}|)}\Big(p_{k|k-1}\big[h(1-p_D+p_DG_Z(0))\big]\Big)\prod_{V\in \mathcal{Q}}\eta_V[0,h]\bigg)\Bigg)\nonumber\\
&\times \sum_{W\in\mathcal{P}}\frac{\frac{\delta}{\delta x}\sum_{\mathcal{Q}\angle W}\zeta_{k|k-1}^{(|\mathcal{Q}|)}\Big(p_{k|k-1}\big[h(1-p_D+p_DG_Z(0))\big]\Big)\prod_{V\in \mathcal{Q}}\eta_V[0,h]}
{\left(\begin{array}{c}\zeta_{FA}^{(|W|)}(0)\\+\sum_{\mathcal{Q}\angle W}\zeta_{k|k-1}^{(|\mathcal{Q}|)}\Big(p_{k|k-1}\big[h(1-p_D+p_DG_Z(0))\big]\Big)\prod_{V\in \mathcal{Q}}\eta_V[0,h]\end{array}\right)}\\
=&\frac{\delta}{\delta Z}F[0,h]\zeta_{k|k-1}^{(1)}\Big(p_{k|k-1}\big[h(1-p_D+p_DG_Z(0))\big]\Big)\nonumber\\
&\times\big(1-p_D(x)+p_D(x)G_Z(0)\big)p_{k|k-1}(x)\nonumber\\
&+F[0,h]\bigg(\prod_{z'\in Z}p_{FA}(z')\bigg)\sum_{\mathcal{P}\angle Z}\Bigg(\prod_{W\in\mathcal{P}}\bigg(\zeta_{FA}^{(|W|)}(0)\nonumber\\
&+\sum_{\mathcal{Q}\angle W}\zeta_{k|k-1}^{(|\mathcal{Q}|)}\Big(p_{k|k-1}\big[h(1-p_D+p_DG_Z(0))\big]\Big)\prod_{V\in \mathcal{Q}}\eta_V[0,h]\bigg)\Bigg)\nonumber\\
&\times \sum_{W\in\mathcal{P}}
\frac{\left(\begin{array}{c}\sum_{\mathcal{Q}\angle W}\Big(\prod_{V\in \mathcal{Q}}\eta_V[0,h]\Big)\bigg(\zeta_{k|k-1}^{(|\mathcal{Q}|+1)}\Big(p_{k|k-1}\big[h(1-p_D+p_DG_Z(0))\big]\Big)\\\times(1-p_D(x)+p_D(x)G_Z(0))p_{k|k-1}(x)\\
+\zeta_{k|k-1}^{(|\mathcal{Q}|)}\Big(p_{k|k-1}\big[h(1-p_D+p_DG_Z(0))\big]\Big)\\\times\sum_{V\in\mathcal{Q}}\frac{p_D(x)G_Z^{(|V|)}(0)\prod_{z'\in V}\frac{p_z(z'|x)}{p_{FA}(z')} p_{k|k-1}(x)}{\eta_{V}[0,h]} \bigg)\end{array}\right)}
{\left(\begin{array}{c}\zeta_{FA}^{(|W|)}(0)\\+\sum_{\mathcal{Q}\angle W}\zeta_{k|k-1}^{(|\mathcal{Q}|)}\Big(p_{k|k-1}\big[h(1-p_D+p_DG_Z(0))\big]\Big)\prod_{V\in \mathcal{Q}}\eta_V[0,h]\end{array}\right)}
\end{align}
Evaluating at $h=1$, we get
\begin{align}
\frac{\delta}{\delta x}\frac{\delta}{\delta Z}F[0,1]=&\frac{\delta}{\delta Z}F[0,1]\zeta_{k|k-1}^{(1)}\Big(p_{k|k-1}\big[1-p_D+p_DG_Z(0)\big]\Big)\big(1-p_D(x)+p_D(x)G_Z(0)\big)p_{k|k-1}(x)\nonumber\\
&+F[0,1]\bigg(\prod_{z'\in Z}p_{FA}(z')\bigg)\sum_{\mathcal{P}\angle Z}\Bigg(\prod_{W\in\mathcal{P}}\bigg(\zeta_{FA}^{(|W|)}(0)\nonumber\\
&+\sum_{\mathcal{Q}\angle W}\zeta_{k|k-1}^{(|\mathcal{Q}|)}\Big(p_{k|k-1}\big[1-p_D+p_DG_Z(0)\big]\Big)\prod_{V\in \mathcal{Q}}\eta_V[0,1]\bigg)\Bigg)\nonumber\\
&\times \sum_{W\in\mathcal{P}}
\frac{\left(\begin{array}{c}\sum_{\mathcal{Q}\angle W}\Big(\prod_{V\in \mathcal{Q}}\eta_V[0,1]\Big)\bigg(\zeta_{k|k-1}^{(|\mathcal{Q}|+1)}\Big(p_{k|k-1}\big[1-p_D+p_DG_Z(0)\big]\Big)\\\times\big(1-p_D(x)+p_D(x)G_Z(0)\big)p_{k|k-1}(x)\\
+\zeta_{k|k-1}^{(|\mathcal{Q}|)}\Big(p_{k|k-1}\big[1-p_D+p_DG_Z(0)\big]\Big)\\\times\sum_{V\in\mathcal{Q}}\frac{p_D(x)G_Z^{(|V|)}(0)\prod_{z'\in V}\frac{p_z(z'|x)}{p_{FA}(z')} p_{k|k-1}(x)}{\eta_{V}[0,1]} \bigg)\end{array}\right)}
{\left(\begin{array}{c}\zeta_{FA}^{(|W|)}(0)\\+\sum_{\mathcal{Q}\angle W}\zeta_{k|k-1}^{(|\mathcal{Q}|)}\Big(p_{k|k-1}\big[1-p_D+p_DG_Z(0)\big]\Big)\prod_{V\in \mathcal{Q}}\eta_V[0,1]\end{array}\right)}
\end{align}
Defining additional quantities
\begin{align}
\alpha_{\mathcal{Q}}\triangleq&\prod_{V\in \mathcal{Q}}\eta_V[0,1]\\
\beta_{W}\triangleq&\zeta_{FA}^{(|W|)}(0)+\sum_{\mathcal{Q}\angle W}\zeta_{k|k-1}^{(|\mathcal{Q}|)}\Big(p_{k|k-1}\big[1-p_D+p_DG_Z(0)\big]\Big)\alpha_\mathcal{Q}
\end{align}
we get
\begin{align}
\frac{\delta}{\delta x}\frac{\delta}{\delta Z}F[0,1]=&\frac{\delta}{\delta Z}F[0,1]\zeta_{k|k-1}^{(1)}\Big(p_{k|k-1}\big[1-p_D+p_DG_Z(0)\big]\Big)\big(1-p_D(x)+p_D(x)G_Z(0)\big)p_{k|k-1}(x)\nonumber\\
&+F[0,1]\bigg(\prod_{z'\in Z}p_{FA}(z')\bigg)\sum_{\mathcal{P}\angle Z}\bigg(\prod_{W\in\mathcal{P}}\beta_W\bigg)\nonumber\\
&\times \sum_{W\in\mathcal{P}}
\frac{1}{\beta_{W}}\left(\begin{array}{c}\sum_{\mathcal{Q}\angle W}\alpha_\mathcal{Q}\bigg(\zeta_{k|k-1}^{(|\mathcal{Q}|+1)}\Big(p_{k|k-1}\big[1-p_D+p_DG_Z(0)\big]\Big)\\\times\big(1-p_D(x)+p_D(x)G_Z(0)\big)p_{k|k-1}(x)\\
+\zeta_{k|k-1}^{(|\mathcal{Q}|)}\Big(p_{k|k-1}\big[1-p_D+p_DG_Z(0)\big]\Big)\\\times\sum_{V\in\mathcal{Q}}\frac{p_D(x)G_Z^{(|V|)}(0)\prod_{z'\in V}\frac{p_z(z'|x)}{p_{FA}(z')} p_{k|k-1}(x)}{\eta_{V}[0,1]} \bigg)\end{array}\right)
\end{align}
Now dividing by $\frac{\delta}{\delta Z}F[0,1]$ to obtain
$D_{k|k}(x)$, we get
\begin{align}
D_{k|k}(x)=&\big(1-p_D(x)+p_D(x)G_Z(0)\big)\zeta_{k|k-1}^{(1)}\Big(p_{k|k-1}\big[1-p_D+p_DG_Z(0)\big]\Big) p_{k|k-1}(x)\nonumber\\
&+\frac{\sum_{\mathcal{P}\angle Z}\bigg(\prod_{W\in\mathcal{P}}\beta_W\bigg)\sum_{W\in\mathcal{P}}
\frac{1}{\beta_{W}}\left(\begin{array}{c}\sum_{\mathcal{Q}\angle W}\alpha_\mathcal{Q}\bigg(\zeta_{k|k-1}^{(|\mathcal{Q}|+1)}\Big(p_{k|k-1}\big[1-p_D+p_DG_Z(0)\big]\Big)\\\times\big(1-p_D(x)+p_D(x)G_Z(0)\big)p_{k|k-1}(x)\\
+\zeta_{k|k-1}^{(|\mathcal{Q}|)}\Big(p_{k|k-1}\big[1-p_D+p_DG_Z(0)\big]\Big)\\\times\sum_{V\in\mathcal{Q}}\frac{p_D(x)G_Z^{(|V|)}(0)\prod_{z'\in V}\frac{p_z(z'|x)}{p_{FA}(z')} p_{k|k-1}(x)}{\eta_{V}[0,1]} \bigg)\end{array}\right)}
{\sum_{\mathcal{P}\angle Z}\prod_{W\in\mathcal{P}}\beta_W}
\end{align}
Defining coefficients $\omega_{\mathcal{P}}$ as
\begin{align}
\omega_{\mathcal{P}}\triangleq\frac{\prod_{W\in\mathcal{P}}\beta_W}{\sum_{\mathcal{P}\angle Z}\prod_{W\in\mathcal{P}}\beta_W}
\end{align}
we get
\begin{align}
D_{k|k}(x)=&\big(1-p_D(x)+p_D(x)G_Z(0)\big)\zeta_{k|k-1}^{(1)}\Big(p_{k|k-1}\big[1-p_D+p_DG_Z(0)\big]\Big) p_{k|k-1}(x)\nonumber\\
&+\sum_{\mathcal{P}\angle Z}\omega_\mathcal{P}\sum_{W\in\mathcal{P}}
\frac{1}{\beta_{W}}\left(\begin{array}{c}\sum_{\mathcal{Q}\angle W}\alpha_\mathcal{Q}\bigg(\zeta_{k|k-1}^{(|\mathcal{Q}|+1)}\Big(p_{k|k-1}\big[1-p_D+p_DG_Z(0)\big]\Big)\\\times\big(1-p_D(x)+p_D(x)G_Z(0)\big)p_{k|k-1}(x)\\
+\zeta_{k|k-1}^{(|\mathcal{Q}|)}\Big(p_{k|k-1}\big[1-p_D+p_DG_Z(0)\big]\Big)\\\times\sum_{V\in\mathcal{Q}}\frac{p_D(x)G_Z^{(|V|)}(0)\prod_{z'\in V}\frac{p_z(z'|x)}{p_{FA}(z')} p_{k|k-1}(x)}{\eta_{V}[0,1]} \bigg)\end{array}\right)\\
=&\big(1-p_D(x)+p_D(x)G_Z(0)\big)\zeta_{k|k-1}^{(1)}\Big(p_{k|k-1}\big[1-p_D+p_DG_Z(0)\big]\Big) p_{k|k-1}(x)\nonumber\\
&+\sum_{\mathcal{P}\angle Z}\omega_\mathcal{P}\sum_{W\in\mathcal{P}}
\frac{1}{\beta_{W}}\left(\begin{array}{c}\sum_{\mathcal{Q}\angle W}\alpha_\mathcal{Q}\bigg(\zeta_{k|k-1}^{(|\mathcal{Q}|+1)}\Big(p_{k|k-1}\big[1-p_D+p_DG_Z(0)\big]\Big)\\\times\big(1-p_D(x)+p_D(x)G_Z(0)\big)\\
+\zeta_{k|k-1}^{(|\mathcal{Q}|)}\Big(p_{k|k-1}\big[1-p_D+p_DG_Z(0)\big]\Big)\\\times\sum_{V\in\mathcal{Q}}\frac{p_D(x)G_Z^{(|V|)}(0)\prod_{z'\in V}\frac{p_z(z'|x)}{p_{FA}(z')}}{\eta_{V}[0,1]} \bigg)\end{array}\right)p_{k|k-1}(x)
\end{align}
If now we define a constant
\begin{align}
\kappa\triangleq\sum_{\mathcal{P}\angle Z}\omega_\mathcal{P}\sum_{W\in\mathcal{P}}\frac{1}{\beta_{W}}\sum_{\mathcal{Q}\angle W}\alpha_\mathcal{Q}\zeta_{k|k-1}^{(|\mathcal{Q}|+1)}\Big(p_{k|k-1}\big[1-p_D+p_DG_Z(0)\big]\Big)
\end{align}
we get
\begin{align}
D_{k|k}(x)=&(\zeta_{k|k-1}^{(1)}+\kappa)\big(1-p_D(x)+p_D(x)G_Z(0)\big)p_{k|k-1}(x)\nonumber\\
&+\sum_{\mathcal{P}\angle Z}\omega_\mathcal{P}\sum_{W\in\mathcal{P}}
\frac{1}{\beta_{W}}\sum_{\mathcal{Q}\angle W}\alpha_\mathcal{Q}\zeta_{k|k-1}^{(|\mathcal{Q}|)}\sum_{V\in\mathcal{Q}}\frac{p_D(x)G_Z^{(|V|)}(0)}{\eta_{V}[0,1]}\prod_{z'\in V}\frac{p_z(z'|x)}{p_{FA}(z')} p_{k|k-1}(x)
\end{align}
where we have dropped the argument of the terms
$\zeta_{k|k-1}^{(\cdot)}\Big(p_{k|k-1}\big[1-p_D+p_DG_Z(0)\big]\Big)$
for simplicity.

We have $G_{k|k}[h]$ as
\begin{align}
G_{k|k}[h]=&\frac{\left(\begin{array}{c}F[0,h]\sum_{\mathcal{P}\angle Z}\prod_{W\in\mathcal{P}}\bigg(\zeta_{FA}^{(|W|)}(0)\\
+\sum_{\mathcal{Q}\angle W}\zeta_{k|k-1}^{(|\mathcal{Q}|)}\Big(p_{k|k-1}\big[h(1-p_D+p_DG_Z(0))\big]\Big)\prod_{V\in \mathcal{Q}}\eta_V[0,h] \bigg)\end{array}\right)}
{\left(\begin{array}{c}F[0,1]\sum_{\mathcal{P}\angle Z}\prod_{W\in\mathcal{P}}\bigg(\zeta_{FA}^{(|W|)}(0)\\
+\sum_{\mathcal{Q}\angle W}\zeta_{k|k-1}^{(|\mathcal{Q}|)}\Big(p_{k|k-1}\big[1-p_D+p_DG_Z(0)\big]\Big)\prod_{V\in \mathcal{Q}}\eta_V[0,1] \bigg)\end{array}\right)}\\
=&\frac{G_{k|k-1}\Big(p_{k|k-1}\big[h(1-p_D+p_DG_Z(0))\big]\Big)}{G_{k|k-1}\Big(p_{k|k-1}\big[1-p_D+p_DG_Z(0)\big]\Big)}\nonumber\\
&\times\frac{\sum_{\mathcal{P}\angle Z}\prod_{W\in\mathcal{P}}\bigg(\zeta_{FA}^{(|W|)}(0)+\sum_{\mathcal{Q}\angle W}\zeta_{k|k-1}^{(|\mathcal{Q}|)}\Big(p_{k|k-1}\big[h(1-p_D+p_DG_Z(0))\big]\Big)\prod_{V\in \mathcal{Q}}\eta_V[0,h] \bigg)}
{\sum_{\mathcal{P}\angle Z}\prod_{W\in\mathcal{P}}\beta_W}
\end{align}
Then,
\begin{align}
G_{k|k}(x)=&\frac{G_{k|k-1}\Big(xp_{k|k-1}\big[1-p_D+p_DG_Z(0)\big]\Big)}{G_{k|k-1}\Big(p_{k|k-1}\big[1-p_D+p_DG_Z(0)\big]\Big)}\nonumber\\
&\times\frac{\sum_{\mathcal{P}\angle Z}\prod_{W\in\mathcal{P}}\bigg(\zeta_{FA}^{(|W|)}(0)+\sum_{\mathcal{Q}\angle W}\alpha_{\mathcal{Q}}x^{|\mathcal{Q}|}\zeta_{k|k-1}^{(|\mathcal{Q}|)}\Big(xp_{k|k-1}\big[1-p_D+p_DG_Z(0)\big]\Big)\bigg)}
{\sum_{\mathcal{P}\angle Z}\prod_{W\in\mathcal{P}}\beta_W}\\
=&\frac{G_{k|k-1}\Big(xp_{k|k-1}\big[\cdot\big]\Big)}{G_{k|k-1}\Big(p_{k|k-1}\big[\cdot\big]\Big)}\frac{\sum_{\mathcal{P}\angle Z}\prod_{W\in\mathcal{P}}\bigg(\zeta_{FA}^{(|W|)}(0)+\sum_{\mathcal{Q}\angle W}\alpha_{\mathcal{Q}}x^{|\mathcal{Q}|}\zeta_{k|k-1}^{(|\mathcal{Q}|)}\Big(xp_{k|k-1}\big[\cdot\big]\Big)\bigg)}
{\sum_{\mathcal{P}\angle Z}\prod_{W\in\mathcal{P}}\beta_W}
\end{align}
where we skipped the arguments of $p_{k|k-1}\big[1-p_D+p_DG_Z(0)\big]$ for the sake of clarity.\newpage
\section{Summary}
Calculate the quantities $D_{k|k}(x)$ and $G_{k|k}(x)$ from the given quantities $D_{k|k-1}(x)$ and $G_{k|k-1}(x)$ as follows.
\begin{align}
D_{k|k}(x)=&(\zeta_{k|k-1}^{(1)}+\kappa)\big(1-p_D(x)+p_D(x)G_Z(0)\big)p_{k|k-1}(x)\nonumber\\
&+\sum_{\mathcal{P}\angle Z}\omega_\mathcal{P}\sum_{W\in\mathcal{P}}
\frac{1}{\beta_{W}}\sum_{\mathcal{Q}\angle W}\alpha_\mathcal{Q}\zeta_{k|k-1}^{(|\mathcal{Q}|)}\sum_{V\in\mathcal{Q}}\frac{p_D(x)G_Z^{(|V|)}(0)}{\eta_{V}[0,1]}\prod_{z'\in V}\frac{p_z(z'|x)}{p_{FA}(z')} p_{k|k-1}(x)\\
G_{k|k}(x)=&\frac{G_{k|k-1}\Big(xp_{k|k-1}\big[\cdot\big]\Big)}{G_{k|k-1}\Big(p_{k|k-1}\big[\cdot\big]\Big)}\frac{\sum_{\mathcal{P}\angle Z}\prod_{W\in\mathcal{P}}\bigg(\zeta_{FA}^{(|W|)}(0)+\sum_{\mathcal{Q}\angle W}\alpha_{\mathcal{Q}}x^{|\mathcal{Q}|}\zeta_{k|k-1}^{(|\mathcal{Q}|)}\Big(xp_{k|k-1}\big[\cdot\big]\Big)\bigg)}
{\sum_{\mathcal{P}\angle Z}\prod_{W\in\mathcal{P}}\beta_W}
\end{align}
where we skipped the arguments of $p_{k|k-1}\big[1-p_D+p_DG_Z(0)\big]$ for the sake of clarity and
\begin{align}
p_{k|k-1}(x)=&\frac{D_{k|k-1}(x)}{\int D_{k|k-1}(x)\d x}\\
\zeta^{(i)}_{FA}(x)\triangleq&\left(\frac{G_{FA}^{(1)}(x)}{G_{FA}(x)}\right)^{(i-1)}\\
\zeta_{k|k-1}^{(i)}(x)\triangleq&\left(\frac{G_{k|k-1}^{(1)}(x)}{G_{k|k-1}(x)}\right)^{(i-1)}\\
\kappa\triangleq&\sum_{\mathcal{P}\angle Z}\omega_\mathcal{P}\sum_{W\in\mathcal{P}}\frac{1}{\beta_{W}}\sum_{\mathcal{Q}\angle W}\alpha_\mathcal{Q}\zeta_{k|k-1}^{(|\mathcal{Q}|+1)}\Big(p_{k|k-1}\big[1-p_D+p_DG_Z(0)\big]\Big)\\
\omega_{\mathcal{P}}\triangleq&\frac{\prod_{W\in\mathcal{P}}\beta_W}{\sum_{\mathcal{P}\angle Z}\prod_{W\in\mathcal{P}}\beta_W}\\
\beta_{W}\triangleq&\zeta_{FA}^{(|W|)}(0)+\sum_{\mathcal{Q}\angle W}\zeta_{k|k-1}^{(|\mathcal{Q}|)}\Big(p_{k|k-1}\big[1-p_D+p_DG_Z(0)\big]\Big)\alpha_\mathcal{Q}\\
\alpha_{\mathcal{Q}}\triangleq&\prod_{V\in \mathcal{Q}}\eta_V[0,1]\\
\eta_V[g,h]\triangleq&p_{k|k-1}\Big[hp_DG_Z^{(|V|)}(p_z[g])\prod_{z'\in V}\frac{p_z(z')}{p_{FA}(z')}\Big]
\end{align}
\newpage
\section{Calculating Posterior Cardinality Distribution $P_{k|k}(n)$}
In this section, we are going to calculate the posterior cardinality
distribution $P_{k|k}(n)$ based on the posterior cardinality
probability generating function (p.g.f.) $G_{k|k}(x)$. For any
cardinality distribution, p.g.f.  pair given as $P(n)$--$G(x)$ we
know the relationship
\begin{align}
P(n)=\frac{1}{n!}G^{(n)}(0)
\end{align}
which we are going to use in this calculation. In the previous
sections, we calculated the posterior cardinality p.g.f.
$G_{k|k}(x)$ as
\begin{align}
G_{k|k}(x)=&\frac{G_{k|k-1}\Big(xp_{k|k-1}\big[\cdot\big]\Big)}{G_{k|k-1}\Big(p_{k|k-1}\big[\cdot\big]\Big)}\frac{\sum_{\mathcal{P}\angle
Z}\prod_{W\in\mathcal{P}}\bigg(\zeta_{FA}^{(|W|)}(0)+\sum_{\mathcal{Q}\angle
W}\alpha_{\mathcal{Q}}x^{|\mathcal{Q}|}\zeta_{k|k-1}^{(|\mathcal{Q}|)}\Big(xp_{k|k-1}\big[\cdot\big]\Big)\bigg)}
{\sum_{\mathcal{P}\angle Z}\prod_{W\in\mathcal{P}}\beta_W}
\end{align}
Here, the denominators on the right hand side are constant with
respect to $x$ and for this reason, we are going to take the
derivative of only the numerator terms.

For the derivative of the multiplication of two functions, we have
\begin{align}
\big(f(x)g(x)\big)^{(n)}=\sum_{i=0}^n C_{n,i}f^{(n-i)}(x)g^{(i)}(x)\label{Eqn:DerivativeTwoMultiplication}
\end{align}
where $C_{n|i}$ is the binomial coefficient given as \begin{align}
C_{n|i}=\frac{n!}{i!(n-i)!}.
\end{align}
We have the generalization of this derivative rule as
\begin{align}
\left(\prod_{j=1}^{M}f_j(x)\right)^{(n)}=\sum_{\substack{0\le i_1,i_2,\ldots,i_M \le M \\ i_1+i_2+\cdots+i_M=n}}C_{n|i_1,\ldots,i_M}\prod_{j=1}^{M}f_j^{(i_j)}(x)\label{Eqn:DerivativeMultipleMultiplication}
\end{align}
where $C_{n|i_1,\ldots,i_M}$ is the multinomial coefficient given as
\begin{align}
C_{n|i_1,\ldots,i_M}=\frac{n!}{i_1!i_2!\cdots i_M!}.
\end{align}
We are going to use the following lemmas in the calculation.
\begin{lemma} Derivatives of the term
$G_{k|k-1}\Big(xp_{k|k-1}\big[\cdot\big]\Big)$ are given as
\begin{align}
G_{k|k-1}^{(n)}\Big(xp_{k|k-1}\big[\cdot\big]\Big)\Big|_{x=0}=&\Big(p_{k|k-1}\big[\cdot\big]\Big)^{n}G_{k|k-1}^{(n)}(0)\\
=&n!\Big(p_{k|k-1}\big[\cdot\big]\Big)^{n}P_{k|k-1}(n)
\end{align}
where $P_{k|k-1}(\cdot)$ is the predicted cardinality
distribution.\hfill$\square$
\end{lemma}
\begin{lemma} Derivatives of the term $\bigg(\zeta_{FA}^{(|W|)}(0)+\sum_{\mathcal{Q}\angle
W}\alpha_{\mathcal{Q}}x^{|\mathcal{Q}|}\zeta_{k|k-1}^{(|\mathcal{Q}|)}\Big(xp_{k|k-1}\big[\cdot\big]\Big)\bigg)$
are given as
\begin{align}
\bigg(\zeta_{FA}^{(|W|)}(0)+\sum_{\mathcal{Q}\angle
W}\alpha_{\mathcal{Q}}x^{|\mathcal{Q}|}\zeta_{k|k-1}^{(|\mathcal{Q}|)}\Big(xp_{k|k-1}\big[\cdot\big]\Big)\bigg)^{(n)}\bigg|_{x=0}=\zeta_{FA}^{(|W|)}(0)\delta_{n,0}+n!\sum_{\substack{\mathcal{Q}\angle W \\ |\mathcal{Q}|=n}}\alpha_\mathcal{Q}\zeta_{k|k-1}^{(|\mathcal{Q}|)}(0)
\end{align}
Notice that the derivative is equal to $0$ if
$n>|W|$.\hfill$\square$
\end{lemma}
\begin{lemma} The $n$th derivative of the summation
\begin{align}
\sum_{\mathcal{P}\angle
Z}\prod_{W\in\mathcal{P}}\bigg(\zeta_{FA}^{(|W|)}(0)+\sum_{\mathcal{Q}\angle
W}\alpha_{\mathcal{Q}}x^{|\mathcal{Q}|}\zeta_{k|k-1}^{(|\mathcal{Q}|)}\Big(xp_{k|k-1}\big[\cdot\big]\Big)\bigg)
\end{align}
is given as
\begin{align}
&\sum_{\substack{\mathcal{P}\angle
Z \\ \mathcal{P}\triangleq\{W_{j}\}_{j=1}^{|\mathcal{P}|}}}\sum_{\substack{0\le i_1,i_2,\ldots,i_{|\mathcal{P}|} \le |\mathcal{P}| \\ i_1+i_2+\cdots+i_{|\mathcal{P}|}=n \\ i_j\le |W_j|\quad\forall j}}C_{n|i_1,\ldots,i_{|\mathcal{P}|}}\prod_{j=1}^{|\mathcal{P}|}
\bigg(\zeta_{FA}^{(|W_j|)}(0)\delta_{i_j,0}+i_j!\sum_{\mathcal{Q}\angle W_j \atop |\mathcal{Q}|=i_j}\alpha_\mathcal{Q}\zeta_{k|k-1}^{(|\mathcal{Q}|)}(0)\bigg)\\
&=n!\sum_{\substack{\mathcal{P}\angle
Z \\ \mathcal{P}\triangleq\{W_{j}\}_{j=1}^{|\mathcal{P}|}}}\sum_{\substack{0\le i_1,i_2,\ldots,i_{|\mathcal{P}|} \le |\mathcal{P}| \\ i_1+i_2+\cdots+i_{|\mathcal{P}|}=n \\ i_j\le |W_j|\quad\forall j}}\prod_{j=1}^{|\mathcal{P}|}
\bigg(\zeta_{FA}^{(|W_j|)}(0)\delta_{i_j,0}+\sum_{\mathcal{Q}\angle W_j \atop |\mathcal{Q}|=i_j}\alpha_\mathcal{Q}\zeta_{k|k-1}^{(|\mathcal{Q}|)}(0)\bigg)
\end{align}
\end{lemma}
where we used \eqref{Eqn:DerivativeMultipleMultiplication} and the
condition $i_j< |W_j|\quad\forall j$ of the second summation is
added just to reduce the number of possibilities.\hfill$\square$

Now we are ready for the final result using the derivative of the
multiplication rule \eqref{Eqn:DerivativeTwoMultiplication} and the
three lemmas above as follows.
\begin{align}
P_{k|k}(n)\triangleq&\frac{1}{G_{k|k-1}\Big(p_{k|k-1}\big[\cdot\big]\Big)\sum_{\mathcal{P}\angle Z}\prod_{W\in\mathcal{P}}\beta_W}\sum_{i=0}^{n}\Big(p_{k|k-1}\big[\cdot\big]\Big)^{n-i}P_{k|k-1}(n-i)\nonumber\\
&\times\sum_{\substack{\mathcal{P}\angle
Z \\ \mathcal{P}\triangleq\{W_{j}\}_{j=1}^{|\mathcal{P}|}}}\sum_{\substack{0\le i_1,i_2,\ldots,i_{|\mathcal{P}|} \le |\mathcal{P}| \\ i_1+i_2+\cdots+i_{|\mathcal{P}|}=i \\ i_j\le |W_j|\quad\forall j}}\prod_{j=1}^{|\mathcal{P}|}
\bigg(\zeta_{FA}^{(|W_j|)}(0)\delta_{i_j,0}+\sum_{\mathcal{Q}\angle W_j \atop |\mathcal{Q}|=i_j}\alpha_\mathcal{Q}\zeta_{k|k-1}^{(|\mathcal{Q}|)}(0)\bigg)
\end{align}

\section{Reducing CPHD for extended targets to PHD for extended targets}
Now we check whether the CPHD equations given above reduce to PHD
when the cluster processes are replaced by Poisson processes. This
is done also for checking the formulas given above. First we give
the related quantities.
\begin{itemize}
\item The extended target measurements are distributed according to
a Poisson process. The corresponding likelihood is given as
\begin{align}
f(Z|x)=n!P_z(n|x)\prod_{z\in Z}p_z(z|x)
\end{align}
where
\begin{align}
P_z(n|x)=\frac{e^{-\gamma(x)}\gamma^n(x)}{n!}
\end{align}
\item The false alarms are distributed according to a Poisson
process also.
\begin{align}
f(Z_{FA})=n!P_{FA}(n)\prod_{z\in Z_{FA}} p_{FA}(z)
\end{align}
where
\begin{align}
P_{FA}(n)=\frac{e^{-\lambda}\lambda^n}{n!}
\end{align}
\item The multitarget prior $f(X_k|Z_{0:k-1})$ is assumed to be a
Poisson process.
\begin{align}
f(X_k|Z_{0:k-1})=n!P_{k+1|k}(n)\prod_{x_k\in X_k} p_{k+1|k}(x_k)
\end{align}
where
\begin{align}
p_{k+1|k}(x_k)\triangleq N_{k|k-1}^{-1}D_{k|k-1}(x_k)
\end{align}
with $N_{k|k-1}\triangleq \int D_{k|k-1}(x_k)\d x_k$. We have
the cardinality distribution $P_{k+1|k}(n)$ given as
\begin{align}
P_{k+1|k}(n)=\frac{e^{-N_{k|k-1}}N_{k|k-1}^n}{n!}
\end{align}
\end{itemize}
We here must note that the probability generating function and
functionals corresponding to a Poisson process with parameter
$\lambda$ and density $p(\cdot)$ corresponds to
\begin{align}
G(x)=\exp(\lambda x-\lambda)\quad\mbox{and}\quad G[h]=\exp(\lambda p[h]-\lambda)
\end{align}
We have
\begin{align}
G^{(n)}(x)=\lambda^nG(x)=\lambda^n\exp(\lambda x-\lambda)
\end{align}
Also the following expression holds.
\begin{align}
G^{(n)}(0)=\lambda^ne^{-\lambda}
\end{align}
The functions $\zeta_{FA}$ and $\zeta_{k|k-1}$ are then given as
\begin{align}
\zeta^{(i)}_{FA}(x)\triangleq&\left(\frac{G_{FA}^{(1)}(x)}{G_{FA}(x)}\right)^{(i-1)}=(\lambda)^{(i-1)}=\begin{cases}0& i>1\\\lambda& i=1\end{cases}=\lambda\delta_{i,1}\\
\zeta_{k|k-1}^{(i)}(x)\triangleq&\left(\frac{G_{k|k-1}^{(1)}(x)}{G_{k|k-1}(x))}\right)^{(i-1)}=(N_{k|k-1})^{(i-1)}=\begin{cases}0& i>1\\ N_{k|k-1}& i=1\end{cases}=N_{k|k-1}\delta_{i,1}
\end{align}
Then in the formula
\begin{align}
\frac{\delta}{\delta Z}F[g,h]=&F[g,h]\bigg(\prod_{z'\in Z}p_{FA}(z')\bigg)\sum_{\mathcal{P}\angle Z}\prod_{W\in\mathcal{P}}\bigg(\zeta_{FA}^{(|W|)}(p_{FA}[g])\nonumber\\
&+\sum_{\mathcal{Q}\angle W}\zeta_{k|k-1}^{(|\mathcal{Q}|)}\Big(p_{k|k-1}\big[h(1-p_D+p_DG_Z(p_z[g]))\big]\Big)\prod_{V\in \mathcal{Q}}\eta_V[g,h]  \bigg),
\end{align}
all terms in the summation $\sum_{\mathcal{Q}\angle W}$ becomes zero
except $\mathcal{Q}=\{W\}$ which means that there is only one set in
the partition $\mathcal{Q}$ and it is $V=W$. Also from above, we have
$\zeta^{(|W|)}_{FA}(x)=\lambda\delta_{|W|,1}$. Substituting these
into the above equation
\begin{align}
\frac{\delta}{\delta Z}F[g,h]=&F[g,h]\bigg(\prod_{z'\in Z}p_{FA}(z')\bigg)\sum_{\mathcal{P}\angle Z}\prod_{W\in\mathcal{P}}\bigg(\lambda\delta_{|W|,1}+N_{k|k-1}\eta_W[g,h]  \bigg)
\end{align}
This equation is the same as the equation Mahler derived for
$\frac{\delta}{\delta Z}F[g,h]$ \cite[Eq. (28)]{mahler_FUSION_2009_extTarg}. In the
CPHD equation
\begin{align}
D_{k|k}(x)=&(\zeta_{k|k-1}^{(1)}+\kappa)\big(1-p_D(x)+p_D(x)G_Z(0)\big)p_{k|k-1}(x)\nonumber\\
&+\sum_{\mathcal{P}\angle Z}\omega_\mathcal{P}\sum_{W\in\mathcal{P}}
\frac{1}{\beta_{W}}\sum_{\mathcal{Q}\angle W}\alpha_\mathcal{Q}\zeta_{k|k-1}^{(|\mathcal{Q}|)}\sum_{V\in\mathcal{Q}}\frac{p_D(x)G_Z^{(|V|)}(0)}{\eta_{V}[0,1]}\prod_{z'\in V}\frac{p_z(z'|x)}{p_{FA}(z')} p_{k|k-1}(x)
\end{align}
By above facts, it is easy to see that $\kappa=0$. Substituting
above facts into this equation we have
\begin{align}
D_{k|k}(x)=&\big(1-p_D(x)+p_D(x)G_Z(0)\big)D_{k|k-1}(x)\nonumber\\
&+\sum_{\mathcal{P}\angle Z}\omega_\mathcal{P}\sum_{W\in\mathcal{P}}
\frac{\alpha_{\{W\}}}{\beta_{W}}\frac{p_D(x)G_Z^{(|W|)}(0)}{\eta_{W}[0,1]}\prod_{z'\in W}\frac{p_z(z'|x)}{p_{FA}(z')} D_{k|k-1}(x)
\end{align}
Seeing that $\alpha_{\{W\}}=\eta_W[0,1]$ we get
\begin{align}
D_{k|k}(x)=&\big(1-p_D(x)+p_D(x)G_Z(0)\big)D_{k|k-1}(x)\nonumber\\
&+\sum_{\mathcal{P}\angle Z}\omega_\mathcal{P}\sum_{W\in\mathcal{P}}
\frac{p_D(x)G_Z^{(|W|)}(0)}{\beta_{W}}\prod_{z'\in W}\frac{p_z(z'|x)}{p_{FA}(z')} D_{k|k-1}(x)
\end{align}
We also see that
\begin{align}
\beta_{W}\triangleq&\zeta_{FA}^{(|W|)}(0)+\sum_{\mathcal{Q}\angle W}\zeta_{k|k-1}^{(|\mathcal{Q}|)}\Big(p_{k|k-1}\big[1-p_D+p_DG_Z(0)\big]\Big)\alpha_\mathcal{Q}\\
=&\lambda\delta_{|W|,1}+N_{k|k-1}\alpha_{\{W\}}\\
=&\lambda\delta_{|W|,1}+N_{k|k-1}\eta_W[0,1]\\
=&\lambda\delta_{|W|,1}+N_{k|k-1}p_{k|k-1}\Big[p_DG_Z^{(|W|)}(0)\prod_{z'\in W}\frac{p_z(z')}{p_{FA}(z')}\Big]\\
=&\lambda\delta_{|W|,1}+D_{k|k-1}\Big[p_DG_Z^{(|W|)}(0)\prod_{z'\in W}\frac{p_z(z')}{p_{FA}(z')}\Big]\\
\end{align}
Dividing and multiplying all quantities in the last summation by $\lambda^{|W|}$ and distributing,
we obtain
\begin{align}
D_{k|k}(x)=&\big(1-p_D(x)+p_D(x)G_Z(0)\big)D_{k|k-1}(x)\nonumber\\
&+\sum_{\mathcal{P}\angle Z}\omega_{\mathcal{P}}\sum_{W\in\mathcal{P}}
\frac{p_D(x)G_Z^{(|W|)}(0)}{\lambda^{-|W|}\beta_{W}}\prod_{z'\in W}\frac{p_z(z'|x)}{\lambda p_{FA}(z')} D_{k|k-1}(x)
\end{align}
Defining the new coefficient
\begin{align}
d_W=&\lambda^{-|W|}\beta_W\\
=&\frac{\lambda}{\lambda^{-|W|}}\delta_{|W|,1}+D_{k|k-1}\Big[p_DG_Z^{(|W|)}(0)\prod_{z'\in W}\frac{p_z(z')}{\lambda p_{FA}(z')}\Big]\\
=&\delta_{|W|,1}+D_{k|k-1}\Big[p_DG_Z^{(|W|)}(0)\prod_{z'\in W}\frac{p_z(z')}{\lambda p_{FA}(z')}\Big]
\end{align}
which are the same coefficients defined in \cite{mahler_FUSION_2009_extTarg}. Also it is easy to see that
\begin{align}
\omega_{\mathcal{P}}\triangleq&\frac{\prod_{W\in\mathcal{P}}\beta_W}{\sum_{\mathcal{P}\angle Z}\prod_{W\in\mathcal{P}}\beta_W}\\
=&\frac{\prod_{W\in\mathcal{P}}\lambda^{-|W|}\beta_W}{\sum_{\mathcal{P}\angle Z}\prod_{W\in\mathcal{P}}\lambda^{-|W|}\beta_W}\\
=&\frac{\prod_{W\in\mathcal{P}}d_W}{\sum_{\mathcal{P}\angle Z}\prod_{W\in\mathcal{P}}d_W}
\end{align}
which are the same $\omega_\mathcal{P}$ coefficients in \cite{mahler_FUSION_2009_extTarg}. Knowing that $G_Z^{|W|}(0)=\gamma^{|W|}(x)e^{-\gamma(x)}$ we get
\begin{align}
D_{k|k}(x)=&\big(1-p_D(x)+p_D(x)e^{-\gamma(x)}\big)D_{k|k-1}(x)\nonumber\\
&+\sum_{\mathcal{P}\angle Z}\omega_{\mathcal{P}}\sum_{W\in\mathcal{P}}
\frac{p_D(x)\gamma^{|W|}(x)e^{-\gamma(x)}}{d_{W}}\prod_{z'\in W}\frac{p_z(z'|x)}{\lambda p_{FA}(z')} D_{k|k-1}(x)
\end{align}
which is the same formula in \cite[Eq. (5)]{mahler_FUSION_2009_extTarg}.
\section{Reducing CPHD for extended targets to CPHD for standard targets}
For standard targets, $P_z(n|x)=0$ except for $n=1$ where it is unity i.e., $P_z(1|x)=1$ which makes $G_Z(x)=x$.
Then
\begin{align}
\eta_V[0,h]\triangleq&p_{k|k-1}\Big[hp_DG_Z^{(|V|)}(0)\prod_{z'\in V}\frac{p_z(z')}{p_{FA}(z')}\Big]\\
=&\begin{cases}p_{k|k-1}\Big[hp_D\frac{p_z(z')}{p_{FA}(z')}\Big]& |V|=1, V={\{z'\}}\\ 0 &\mbox{otherwise}\end{cases}\\
=&p_{k|k-1}\Big[hp_D\frac{p_z(z')}{p_{FA}(z')}\Big]\delta_{|V|,1}
\end{align}
This gives
\begin{align}
\alpha_{\mathcal{Q}}\triangleq&\prod_{V\in \mathcal{Q}}\eta_V[0,1]=\prod_{z\in W}p_{k|k-1}\Big[p_D\frac{p_z(z')}{p_{FA}(z')}\Big]\delta_{|\mathcal{Q}|,|W|}
\end{align}
where $\mathcal{Q}$ is a partition of the set $W$.
\begin{align}
\beta_{W}\triangleq&\zeta_{FA}^{(|W|)}(0)+\sum_{\mathcal{Q}\angle W}\zeta_{k|k-1}^{(|\mathcal{Q}|)}\Big(p_{k|k-1}\big[1-p_D+p_DG_Z(0)\big]\Big)\alpha_\mathcal{Q}\\
=&\zeta_{FA}^{(|W|)}(0)+\zeta_{k|k-1}^{(|W|)}\Big(p_{k|k-1}\big[1-p_D\big]\Big)\prod_{z'\in W}p_{k|k-1}\Big[p_D\frac{p_z(z')}{p_{FA}(z')}\Big]
\end{align}
Substituting these into $G_{k|k}[h]$ we get
\begin{align}
G_{k|k}[h]\triangleq&\frac{G_{k|k-1}\Big(p_{k|k-1}\big[h(1-p_D+p_DG_Z(0))\big]\Big)}{G_{k|k-1}\Big(p_{k|k-1}\big[1-p_D+p_DG_Z(0)\big]\Big)}\nonumber\\
&\times\frac{\sum_{\mathcal{P}\angle Z}\prod_{W\in\mathcal{P}}\bigg(\zeta_{FA}^{(|W|)}(0)+\sum_{\mathcal{Q}\angle W}\zeta_{k|k-1}^{(|\mathcal{Q}|)}\Big(p_{k|k-1}\big[h(1-p_D+p_DG_Z(0))\big]\Big)\prod_{V\in \mathcal{Q}}\eta_V[0,h] \bigg)}
{\sum_{\mathcal{P}\angle Z}\prod_{W\in\mathcal{P}}\beta_W}\\
=&\frac{G_{k|k-1}\Big(p_{k|k-1}\big[h(1-p_D)\big]\Big)}{G_{k|k-1}\Big(p_{k|k-1}\big[1-p_D\big]\Big)}\nonumber\\
&\times\frac{\sum_{\mathcal{P}\angle Z}\prod_{W\in\mathcal{P}}\bigg(\zeta_{FA}^{(|W|)}(0)+\zeta_{k|k-1}^{(|W|)}\Big(p_{k|k-1}\big[h(1-p_D)\big]\Big)\prod_{z'\in W}p_{k|k-1}\Big[hp_D\frac{p_z(z')}{p_{FA}(z')}\Big]\bigg)}
{\sum_{\mathcal{P}\angle Z}\prod_{W\in\mathcal{P}}\bigg(\zeta_{FA}^{(|W|)}(0)+\zeta_{k|k-1}^{(|W|)}\Big(p_{k|k-1}\big[1-p_D\big]\Big)\prod_{z'\in W}p_{k|k-1}\Big[p_D\frac{p_z(z')}{p_{FA}(z')}\Big]\bigg)}\\
=&\frac{G_{k|k-1}\Big(p_{k|k-1}\big[h(1-p_D)\big]\Big)}{G_{k|k-1}\Big(p_{k|k-1}\big[1-p_D\big]\Big)}\nonumber\\
&\times\frac{\sum_{\mathcal{P}\angle Z}\prod_{W\in\mathcal{P}}\frac{\delta}{\delta W}\bigg(\log G_{FA}(p_{FA}[g])+\log G_{k|k-1}\Big(p_{k|k-1}\big[h(1-p_D+p_Dp_z[g])\big]\Big)\bigg)\bigg|_{g=0}}
{\sum_{\mathcal{P}\angle Z}\prod_{W\in\mathcal{P}}\frac{\delta}{\delta W}\bigg(\log G_{FA}(p_{FA}[g])+\log G_{k|k-1}\Big(p_{k|k-1}\big[1-p_D+p_Dp_z[g]\big]\Big)\bigg)\bigg|_{g=0}}\label{Eqn:ReducingECPHDtoCPHDold}\\
=&\frac{G_{k|k-1}\Big(p_{k|k-1}\big[h(1-p_D)\big]\Big)}{G_{k|k-1}\Big(p_{k|k-1}\big[1-p_D\big]\Big)}\nonumber\\
&\times\frac{\sum_{\mathcal{P}\angle Z}\prod_{W\in\mathcal{P}}\frac{\delta}{\delta W}\log\bigg(G_{FA}(p_{FA}[g])G_{k|k-1}\Big(p_{k|k-1}\big[h(1-p_D+p_Dp_z[g])\big]\Big)\bigg)\bigg|_{g=0}}
{\sum_{\mathcal{P}\angle Z}\prod_{W\in\mathcal{P}}\frac{\delta}{\delta W}\log \bigg(G_{FA}(p_{FA}[g])G_{k|k-1}\Big(p_{k|k-1}\big[1-p_D+p_Dp_z[g]\big]\Big)\bigg)\bigg|_{g=0}}\label{Eqn:ReducingECPHDtoCPHD}
\end{align}
where we used the facts that $\int \zeta_{FA}^{(1)}=\log G_{FA}$ and $\int \zeta_{k|k-1}^{(1)}=\log G_{k|k-1}$ to obtain \eqref{Eqn:ReducingECPHDtoCPHDold}.
We can see that for any p.g.fl., $G[g]$
\begin{align}
\sum_{\mathcal{P}\angle Z}\prod_{W\in\mathcal{P}}\left(G^{|W|}[g]+\frac{\delta G}{\delta W}[g]\right)=\frac{\delta G^{|Z|}}{\delta Z}[g]
\end{align}
from the product rule for functional derivatives (See \cite[$8$th row of Table 11.2]{mahler_book_2007}.).
Since $G[g]$ is arbitrary, we can replace it with $G[g]-G[0]$ which gives
\begin{align}
\sum_{\mathcal{P}\angle Z}\prod_{W\in\mathcal{P}}\left((G[g]-G[0])^{|W|}+\frac{\delta G}{\delta W}[g]\right)=\frac{\delta (G[g]-G[0])^{|Z|}}{\delta Z}
\end{align}
Now, evaluating both sides at $g=0$, we get
\begin{align}
\sum_{\mathcal{P}\angle Z}\prod_{W\in\mathcal{P}}\frac{\delta G}{\delta W}[0]=\frac{\delta (G[g]-G[0])^{|Z|}}{\delta Z}\bigg|_{g=0}
\end{align}
Now suppose we have two functionals and we need to evaluate the following ration.
\begin{align}
\frac{\sum_{\mathcal{P}\angle Z}\prod_{W\in\mathcal{P}}\frac{\delta G_1}{\delta W}[0]}{\sum_{\mathcal{P}\angle Z}\prod_{W\in\mathcal{P}}\frac{\delta G_2}{\delta W}[0]}=\frac{\frac{\delta (G_1[g]-G_1[0])^{|Z|}}{\delta Z}\big|_{g=0}}{\frac{\delta (G_2[g]-G_2[0])^{|Z|}}{\delta Z}\big|_{g=0}}
\end{align}
This division gives $\frac{0}{0}$ indeterminate form. One must use l'Hopital's rule $|Z|-1$ times (take derivatives of the numerator and denominator $|Z|-1$ times with respect to $g$),
to get
\begin{align}
\frac{\sum_{\mathcal{P}\angle Z}\prod_{W\in\mathcal{P}}\frac{\delta G_1}{\delta W}[0]}{\sum_{\mathcal{P}\angle Z}\prod_{W\in\mathcal{P}}\frac{\delta G_2}{\delta W}[0]}=\frac{\frac{\delta G_1[g]}{\delta Z}\big|_{g=0}}{\frac{\delta G_2[g]}{\delta Z}\big|_{g=0}}
\end{align}

Using this formula for evaluating \eqref{Eqn:ReducingECPHDtoCPHD}, we obtain
\begin{align}
G_{k|k}[h]=&\frac{G_{k|k-1}\Big(p_{k|k-1}\big[h(1-p_D)\big]\Big)}{G_{k|k-1}\Big(p_{k|k-1}\big[1-p_D\big]\Big)}\nonumber\\
&\times\frac{\frac{\delta}{\delta Z}\log\bigg(G_{FA}(p_{FA}[g])G_{k|k-1}\Big(p_{k|k-1}\big[h(1-p_D+p_Dp_z[g])\big]\Big)\bigg)\bigg|_{g=0}}{\frac{\delta}{\delta Z}\log \bigg(G_{FA}(p_{FA}[g])G_{k|k-1}\Big(p_{k|k-1}\big[1-p_D+p_Dp_z[g]\big]\Big)\bigg)\bigg)\Big|_{g=0}}\\
=&\frac{\frac{\delta}{\delta Z}\left(G_{FA}[g]G_{k|k-1}\big[h(1-p_D+p_Dp_z[g])\big]\right)\big|_{g=0}}
{\frac{\delta}{\delta Z}\left(G_{FA}[g]G_{k|k-1}\big[1-p_D+p_Dp_z[g]\big]\right)\big|_{g=0}}
\end{align}
which is the formula for CPHD for standard targets (See \cite[Equations (113) and (114)]{Mahler:2007}.).
\section{Proof of the Main Equation \eqref{Eqn:FormulaForF}}
\label{Sec:ProofOfDerivative}
The proof is by induction as in \cite{mahler_FUSION_2009_extTarg}. The formula can be seen to be satisfied for
first and second order derivatives in \eqref{Eqn:FirstDerivativeF} and \eqref{Eqn:SecondDerivativeF} respectively.
We assume that for $|Z|=m$ the formula is satisfied as below.
\begin{align}
\frac{\delta}{\delta Z}F[g,h]=&F[g,h]\bigg(\prod_{z'\in Z}p_{FA}(z')\bigg)\sum_{\mathcal{P}\angle Z}\prod_{W\in\mathcal{P}}\bigg(\zeta_{FA}^{(|W|)}(p_{FA}[g])+\sum_{\mathcal{Q}\angle W}\zeta_{k|k-1}^{(|\mathcal{Q}|)}\prod_{V\in \mathcal{Q}}\eta_V[g,h]  \bigg)
\end{align}
Now we need to show that the formula is satisfied for $Z\cup z_\star$.
\begin{align}
\frac{\delta}{\delta z_\star}\frac{\delta}{\delta Z}F[g,h]=&\frac{\delta}{\delta z_\star}\Bigg(F[g,h]\bigg(\prod_{z'\in Z}p_{FA}(z')\bigg)\sum_{\mathcal{P}\angle Z}\prod_{W\in\mathcal{P}}\bigg(\zeta_{FA}^{(|W|)}(p_{FA}[g])+\sum_{\mathcal{Q}\angle W}\zeta_{k|k-1}^{(|\mathcal{Q}|)}\prod_{V\in \mathcal{Q}}\eta_V[g,h]  \bigg)\Bigg)\\
=&F[g,h]\bigg(\prod_{z'\in Z\cup z_\star}p_{FA}(z')\bigg)\bigg(\zeta_{FA}^{(1)}(p_{FA}[g])\nonumber\\
&+\zeta_{k|k-1}^{(1)}\eta_{\{z_\star\}}[g,h]\bigg)\nonumber\\
&\times\sum_{\mathcal{P}\angle Z}\prod_{W\in\mathcal{P}}\bigg(\zeta_{FA}^{(|W|)}(p_{FA}[g])+\sum_{\mathcal{Q}\angle W}\zeta_{k|k-1}^{(|\mathcal{Q}|)}\prod_{V\in \mathcal{Q}}\eta_V[g,h]  \bigg)\nonumber\\
&+F[g,h]\bigg(\prod_{z'\in Z\cup z_\star}p_{FA}(z')\bigg)\sum_{\mathcal{P}\angle Z}\Bigg(\prod_{W\in\mathcal{P}}\bigg(\zeta_{FA}^{(|W|)}(p_{FA}[g])+\sum_{\mathcal{Q}\angle W}\zeta_{k|k-1}^{(|\mathcal{Q}|)}\prod_{V\in \mathcal{Q}}\eta_V[g,h]\bigg)\Bigg)\nonumber\\
&\times\sum_{W\in\mathcal{P}}
\frac{\left(\begin{array}{c}\zeta_{FA}^{(|W|+1)}(p_{FA}[g])\\+\sum_{\mathcal{Q}\angle W}\Big(\prod_{V\in \mathcal{Q}}\eta_V[g,h]\Big)\Big(\zeta_{k|k-1}^{(|\mathcal{Q}|+1)}\eta_{\{z_\star\}}[g,h]+\sum_{V\in\mathcal{Q}}\frac{\eta_{V\cup z_\star}[g,h]}{\eta_V[g,h]}\Big)\bigg)\end{array}\right)}
{\zeta_{FA}^{(|W|)}(p_{FA}[g])+\sum_{\mathcal{Q}\angle W}\zeta_{k|k-1}^{(|\mathcal{Q}|)}\prod_{V\in \mathcal{Q}}\eta_V[g,h]}\\
=&F[g,h]\bigg(\prod_{z'\in Z\cup z_\star}p_{FA}(z')\bigg)\sum_{\mathcal{P}\angle Z}\prod_{W\in\mathcal{P}\cup \{z_\star\}}\bigg(\zeta_{FA}^{(|W|)}(p_{FA}[g])+\sum_{\mathcal{Q}\angle W}\zeta_{k|k-1}^{(|\mathcal{Q}|)}\prod_{V\in \mathcal{Q}}\eta_V[g,h]  \bigg)\nonumber\\
&+F[g,h]\bigg(\prod_{z'\in Z\cup z_\star}p_{FA}(z')\bigg)\sum_{\mathcal{P}\angle Z}\Bigg(\prod_{W\in\mathcal{P}}\bigg(\zeta_{FA}^{(|W|)}(p_{FA}[g])+\sum_{\mathcal{Q}\angle W}\zeta_{k|k-1}^{(|\mathcal{Q}|)}\prod_{V\in \mathcal{Q}}\eta_V[g,h]\bigg)\Bigg)\nonumber\\
&\times\sum_{W\in\mathcal{P}}
\frac{\left(\begin{array}{c}\zeta_{FA}^{(|W|+1)}(p_{FA}[g])+\sum_{\mathcal{Q}\angle W}\zeta_{k|k-1}^{(|\mathcal{Q}|+1)}\prod_{V\in \mathcal{Q}\cup\{z_\star\}}\eta_V[g,h]\\+\sum_{\mathcal{Q}\angle W}\Big(\prod_{V\in \mathcal{Q}}\eta_V[g,h]\Big)\sum_{V\in\mathcal{Q}}\frac{\eta_{V\cup z_\star}[g,h]}{\eta_V[g,h]}\Big)\bigg)\end{array}\right)}
{\zeta_{FA}^{(|W|)}(p_{FA}[g])+\sum_{\mathcal{Q}\angle W}\zeta_{k|k-1}^{(|\mathcal{Q}|)}\prod_{V\in \mathcal{Q}}\eta_V[g,h]}\\
=&F[g,h]\bigg(\prod_{z'\in Z\cup z_\star}p_{FA}(z')\bigg)\sum_{\mathcal{P}\angle Z}\prod_{W\in\mathcal{P}\cup \{z_\star\}}\bigg(\zeta_{FA}^{(|W|)}(p_{FA}[g])+\sum_{\mathcal{Q}\angle W}\zeta_{k|k-1}^{(|\mathcal{Q}|)}\prod_{V\in \mathcal{Q}}\eta_V[g,h]  \bigg)\nonumber\\
&+F[g,h]\bigg(\prod_{z'\in Z\cup z_\star}p_{FA}(z')\bigg)\sum_{\mathcal{P}\angle Z}\Bigg(\prod_{W\in\mathcal{P}}\bigg(\zeta_{FA}^{(|W|)}(p_{FA}[g])+\sum_{\mathcal{Q}\angle W}\zeta_{k|k-1}^{(|\mathcal{Q}|)}\prod_{V\in \mathcal{Q}}\eta_V[g,h]\bigg)\Bigg)\nonumber\\
&\times\sum_{W\in\mathcal{P}}
\frac{\zeta_{FA}^{(|W|+1)}(p_{FA}[g])+\sum_{\mathcal{Q}\angle W\cup\{z_\star\}}\zeta_{k|k-1}^{(|\mathcal{Q}|)}\prod_{V\in \mathcal{Q}}\eta_V[g,h]}
{\zeta_{FA}^{(|W|)}(p_{FA}[g])+\sum_{\mathcal{Q}\angle W}\zeta_{k|k-1}^{(|\mathcal{Q}|)}\prod_{V\in \mathcal{Q}}\eta_V[g,h]}\\
=&F[g,h]\bigg(\prod_{z'\in Z\cup z_\star}p_{FA}(z')\bigg)\sum_{\mathcal{P}\angle Z\cup\{z_\star\}}\prod_{W\in\mathcal{P}}\bigg(\zeta_{FA}^{(|W|)}(p_{FA}[g])+\sum_{\mathcal{Q}\angle W}\zeta_{k|k-1}^{(|\mathcal{Q}|)}\prod_{V\in \mathcal{Q}}\eta_V[g,h]  \bigg)\nonumber\\
\end{align}
which completes the proof.

\bibliographystyle{ieeetran}
\bibliography{IEEEabrv,CPHDextended}
\end{document}